\documentclass{article}
\usepackage[english]{babel}
\usepackage[english]{babel}
\usepackage{lineno} 
\usepackage{graphicx,multicol}
\usepackage{epic,eepic,epsfig}
\usepackage{amssymb}
\usepackage{amsmath,amsthm}
\usepackage{color}
\usepackage{tikz}

\usepackage{eso-pic} 

\newcommand{\smallqed}{\hfill {\tiny ($\Box$)}}

\newcommand{\eqSp}{\vspace{0.1cm}}

\newcommand{\2}{\vspace{2mm}}

\theoremstyle{plain}
\newtheorem{theorem}{Theorem}

\newtheorem{lemma}[theorem]{Lemma}

\theoremstyle{definition}

\newtheorem{conjecture}[theorem]{Conjecture}

\newcommand{\greencolor}{green}
\newcommand{\greenline}{solid}
\newcommand{\redcolor}{red}
\newcommand{\redline}{solid}
\newcommand{\bluecolor}{blue}
\newcommand{\blueline}{solid}
\newcommand{\thickline}{thin}

\begin{document}
\bibliographystyle{plain}

\title{Spanning eulerian subdigraphs avoiding $k$ prescribed arcs in tournaments}
\author{J\o{}rgen Bang-Jensen\thanks{Department of Mathematics and Computer Science, University of southern Denmark (email: jbj@imada.sdu.dk). The research was supported by the danish national research foundation under grant number 7014-00037B}, Hugues D\'epr\'es\thanks{Department of Computer Science, ENS Lyon (email: hugues.depres@ens-lyon.fr)}, Anders Yeo\thanks{Department of Mathematics and Computer Science, University of southern Denmark (email: yeo@imada.sdu.dk). The research was supported by the danish national research foundation under grant number 7014-00037B}}

\maketitle

\begin{abstract}
  A digraph is {\bf eulerian} if it is connected and every vertex has its in-degree equal to its out-degree. 
  Having a spanning eulerian subdigraph is thus a weakening of having a hamiltonian cycle.
  A digraph is {\bf semicomplete} if it has no pair of non-adjacent vertices. A {\bf tournament} is a semicomplete digraph without directed cycles of length 2.
 Fraise and Thomassen \cite{fraisseGC3} proved that every $(k+1)$-strong tournament has a hamiltonian cycle which avoids any prescribed set of $k$ arcs. In \cite{bangsupereuler} the authors demonstrated that a number of results concerning vertex-connectivity and hamiltonian cycles in tournaments and  have analogues when we replace vertex connectivity by arc-connectivity and hamiltonian cycles by spanning eulerian subdigraphs. They  showed the existence of a smallest function $f(k)$ such that every $f(k)$-arc-strong semicomplete digraph has a spanning eulerian subdigraph which avoids any prescribed set of $k$ arcs. They proved that $f(k)\leq \frac{(k+1)^2}{4}+1$ and also proved that $f(k)=k+1$ when $k=2,3$. Based on this they conjectured that $f(k)=k+1$ for all $k\geq 0$. In this paper we prove that $f(k)\leq (\lceil\frac{6k+1}{5}\rceil)$.\\

 \noindent{}{\bf Keywords:} Arc-connectivity, Eulerian subdigraph, Tournament, Semicomplete digraph, avoiding prescribed arcs
\end{abstract}

\section{Introduction}
The terminology is consistent with \cite{bang2009}.
 A  classical result on digraphs is  Camion's Theorem
(it was  originally formulated only for tournaments but   easily extends to semicomplete digraphs). 

\begin{theorem}[Camion~\cite{Cami59}]\label{thm:camion}
Every strongly connected  semicomplete digraph has a hamiltonian cycle.
\end{theorem}

A digraph is {\bf connected} if its
underlying undirected graph is connected.
A digraph $D$ is {\bf eulerian} if it contains a spanning eulerian trail $W$ such that $A(W)=A(D)$, or, equivalently by Euler's theorem, if $D$ is connected and $d^+(v)=d^-(v)$ for all $v\in V(D)$.
Finally we say that $D=(V,A)$ is {\bf supereulerian} if it has a spanning eulerian subdigraph $D'=(V,A')$.
By Camion's theorem, every strong semicomplete digraph is supereulerian. Bang-Jensen and Thomass\'e made the following conjecture in 2011 (see e.g. \cite{bangJGT79}) which may be seen as a generalization of Camion's theorem.

\begin{conjecture}
  \label{conj:BJTh}
  Every digraph $D$ with $\lambda{}(D)\geq \alpha{}(D)$ is supereulerian.
\end{conjecture}

This conjecture is open even for digraphs of independence number 2 but has been verified for a couple of classes of digraphs (\cite{alsatamiAM7,bangJGT79}.

An {\bf eulerian factor} of a digraph $D=(V,A)$ is a spanning subdigraph $H=(V,A')$ so that $d_H^+(v)=d_H^-(v)>0$ for all $v\in V$. 
By a {\bf component} of the eulerian factor $H$  we mean a connected component of the digraph $H$.
Bang-Jensen and Maddaloni proved the following result in support of Conjecture \ref{conj:BJTh}.

\begin{theorem}\cite{bangJGT79}
  \label{lem:efactoravoid}
  Every digraph $D$ with $\lambda{}(D)\geq\alpha{}(D)$ has an eulerian factor. Furthermore, such a factor can be found in polynomial time.
\end{theorem}

  There are many results on hamiltonian cycles in tournaments and semicomplete digraphs
  (see e.g. \cite{bangTchapter,kuehnEJC33}). One of these is the following, due to Fraisse and Thomassen. A digraph $D$ is {\bf $\mathbf{k}$-strong} for some natural number $k\geq 1$ if it has at least $k+1$ vertices and remains strongly connected  whenever we delete at most $k-1$ vertices. 

\begin{theorem}\cite{fraisseGC3}
 \label{thm:FrTh}
  Let $T=(V,A)$ be a $(k+1)$-strong tournament and let $\hat{A}\subset A$ have size $k$. Then $T\setminus \hat{A}$ has a hamiltonian cycle.
  \end{theorem}

  This is best possible in terms of the connectivity since there are infinitely many $(k+1)$-strong tournaments that have a vertex of out-degree $k+1$.



  The authors of \cite{bangsupereuler} showed that a number of results concerning hamiltonian paths and cycles in tournaments and semicomplete digraphs and  vertex-connectivity have analogues in results on spanning eulerian subdigraphs in semicomplete digraphs where we can replace vertex-connectivity by  arc-connectivity. A digraph $D$ is {\bf $\mathbf{k}$-arc-strong} for some natural number $k\geq 1$ if it  is strongly connected and remains so when we delete an arbitrary set of at most $k-1$ of its arcs.

  It was shown in \cite{bangsupereuler} that there exists a smallest function $f(k)$ such that every
  $f(k)$-arc-strong semicomplete digraph has a spanning eulerian subdigraph avoiding any set of $k$ prescribed arcs. They proved that $f(k)\leq \frac{(k+1)^2}{4}+1$ and that $f(2)=3, f(3)=4$.
Note that, if true, Conjecture \ref{conj:BJTh} would imply that $f(k)\leq k+2\sqrt{k}$ since deleting $k$ edges from a complete graph cannot generate an independent set of size more than $2\sqrt{k}$ (the exact bound is not important here).

\begin{lemma}\cite{bangsupereuler}
  \label{lem:factoravoidk}
  Let $D=(V,A)$ be a $(k+1)$-arc-strong semicomplete digraph and let $\hat{A}$ be a set of $k$ arcs from $D$. Then $D\setminus \hat{A}$ has an eulerian factor.\
  \end{lemma}
 
  Based on these observations the authors of \cite{bangsupereuler} posed the following conjecture that can be seen as such an arc-analogue of Theorem \ref{thm:FrTh}.
\begin{conjecture}
    \label{conj:euleravoid}
    Let $D=(V,A)$ be a $(k+1)$-arc-strong semicomplete digraph and let $A'\subset A$ be any set of $k$ arcs of $D$. Then $D\setminus A'$ has is supereulerian.
  \end{conjecture}

In this paper we prove the following theorem thus providing additional support for Conjecture \ref{conj:BJTh}.
  \begin{theorem}
    \label{thm:f(k)}
    Let $D=(V,A)$ be a semicomplete digraph with $\lambda{}(D)\geq \frac{6k+1}{5}$ and let $A\subset A$ be an arbitrary set of $k$ arcs of $D$. Then $D\setminus A'$ is supereulerian.
    \end{theorem}
    
We shall use the following theorem due to Meyniel.

\begin{theorem}\cite{meynielJCT14}
  \label{meyniel}
Let $G$ be a strongly connected digraph on $n \geq 2$ vertices.
If $d^+(x)+d^-(x)+d^+(y)+d^-(y) \geq 2|V(G)| -1$ for all pairs of non-adjacent vertices in $G$, then $G$ has a Hamilton cycle.
\end{theorem}

\section{Proof of Theorem \ref{thm:f(k)}}

%
%

\begin{proof}
Let $D$ and $A'$ be defined as in the theorem and let $D'=D-A$. Note that $D'$ is 
$\frac{k+1}{5}$-arc-strong and let $q  = \frac{k+1}{5}$. By Lemma~\ref{lem:factoravoidk} we note that $D'$ contains an eulerian factor, $F$.
Assume that $F = C_1 \cup C_2 \cup \cdots \cup C_p$, where $C_i$ denotes the components in $F$ and that $p$ is minimum possible.
We may assume that $p \geq 2$ as otherwise we are done.
Any vertex $u \in V(D')$ which is adjacent to all vertices in $V(D')\setminus \{u\}$ is called {\bf universal}. 
We now prove the following claim.

\2

{\bf Claim 1:} {\em There is a universal vertex, $u$, in $D'$}.

\2

{\em Proof of Claim 1:} For the sake of contradiction, assume that there is no universal vertex in $D'$.
By Theorem~\ref{meyniel} there exists two non-adjacent vertices $x$ and $y$ in $D'$ such that
$d_{D'}^+(x)+d_{D'}^-(x)+d_{D'}^+(y)+d_{D'}^-(y) < 2|V(D')|-1$, as otherwise there is a Hamilton cycle in $D'$ and we are done.

Let $A' = \{a_1,a_2,\ldots,a_k\}$. Start with $D$ and remove $a_1$, then remove $a_2$, then $a_3$ and continue
this process until we obtain $D'$.  When removing $e_i$ above let $r_i$ denote how many more vertices become non-universal
(that is, they previously had no non-neighbour, but after the removal of $e_i$ they get a non-neighbour) for $i \in [k]$.
Note that $r_i \in \{0,1,2\}$ as the only vertices that can become non-universal are the endpoints of $e_i$.
Let $n_j$ be the number of $r_i$'s of value j for $j \in \{0,1,2\}$.

As there is no universal vertex in $D'$ the following is the number of non-universal vertices in $D'$.

\[
|V(D')| = 0 n_0 + 1 n_1 + 2 n_2 = n_1 + 2n_2  
\]

Let $f(x,y,D) = d_D^+(x)+d_D^-(x)+d_D^+(y)+d_D^-(y)$ and $f(x,y,D') = d_{D'}^+(x)+d_{D'}^-(x)+d_{D'}^+(y)+d_{D'}^-(y)$.
Note that $f(x,y,D') < 2|V(D')|-1$ and $f(x,y,D) \geq 4 \frac{6k+1}{5}$ (as $D$ is $\frac{6k+1}{5}$-arc-strong, implying
that every vertex has at least $\frac{6k+1}{5}$ arcs into it and out of it).
Therefore, by removing the arcs of $A'$ the $f$-value needs to drop by at least $\frac{24k+4}{5} - (2|V(D')|-1)$.

Without loss of generality we may assume that $a_1$ is an arc between $x$ and $y$, and if there is two arcs
between $x$ and $y$ then $a_2$ is the second such arc. Note that after removing the arcs 
between $x$ and $y$ all other arcs, $a_i$, touching $\{x,y\}$ will have $r_i \in \{0,1\}$.

If there are two arcs between $x$ and $y$ we note that the $r_1=0$ and $r_2=2$ and the
value of $f(x,y,D)$ drops by at most $3+n_0+n_1\leq 2 + 2n_0 + n_1$ (as after removing the first two arcs it has dropped by
$4$ and we have one arc with $r$-value $0$).

If there is only one arc between $x$ and $y$ we note that $r_1=2$  and the
value of $f(x,y,D)$ drops by at most $2 + n_0 + n_1$ (as after removing the first arc it has dropped by
$2$ and we have one arc with $r$-value $2$).
Therefore the following holds.

\[
2 + 2n_0 + n_1 \geq f(x,y,D) - f(x,y,D') = \frac{24k+4}{5} - (2|V(D')|-1)
\]

We have previously shown that $|V(D')| = n_1 + 2n_2$, which by the above implies the following.

\[
2 + 2n_0 + n_1 \geq \frac{24k+4}{5} - (2(n_1+2n_2)-1)
\]

Simplifying the above gives us the following, as $n_0+n_1+n_2 = k$.

\[
4k = 4(n_0+n_1+n_2) \geq  2n_0 + 3n_1 + 4n_2  \geq \frac{24k-1}{5} 
\]

However the above is a contradiction, as $20 k \geq 24k-1$ is impossible when $k \geq 1$.~\smallqed{}

\2

By Claim~1, we will assume that $F$ is chosen such that $C_1$ contains the universal vertex, $u$, given by Claim~1 and
that $|V(C_p)|$ is as large as possible. That is, we want $|V(C_p)|$ to be as large as possible, without containing all
universal vertices.
We call a vertex in $V\setminus V(C_p)$ {\bf $C_p$-out-universal} if it dominates all vertices in $C_p$ and analogously we call
it {\bf $C_p$-in-universal} if it is dominated by all vertices in $C_p$.
We call a vertex $u \in V(L)$ {\bf $C_p$-almost-out-universal} if it dominates all vertices in $C_p$ except at most  one vertex.
We now have the following observation.

\2

{\bf Observation 2:} {\em Any vertex $w \in V(D') \setminus V(C_p)$ with no non-neighbour in $C_p$ (in $D'$) is
either $C_p$-out-universal or $C_p$-in-universal.

We may assume without loss of generality that $u$ (defined in Claim~1) is $C_p$-out-universal.}

\2

{\em Proof of Observation 2:} Let $w \in V(D') \setminus V(C_p)$ have no non-neighbour in $C_p$ (in $D'$),
but for the sake of contradiction assume that $w$ is not $C_p$-out-universal and not $C_p$-in-universal.
As $w$ is not $C_p$-out-universal there exists a $x \in V(C_p)$, such that  $xw \in A(D')$.
Starting in $x$ and moving along an eulerian trail in $C_p$ we continue until we come to a vertex $v$ which 
does not dominate $w$ (which exists since $w$ is not $C_p$-in-universal). So there exists an arc $x'v \in A(C_p)$ such that
$x'w \in A(D')$ but $vw \not\in A(D')$. However as $w$ is universal we must then have that $wv \in A(D')$.
Now removing the arc $x'v$ and adding the arcs $x'w$ and $wv$ gives us an eulerian factor in $D'$ with fewer 
components, a contradiction. This completes the first part of the proof of Observation 2.

Without loss of generality we may assume that $u$ is $C_p$-out-universal, as if it is is $C_p$-in-universal
we can reverse all arcs.~\smallqed{}

\2

Let $L_i$ be an inbranching in $C_i$ containing only $F$-arcs for all $i \in [p-1]$, such that the root of $L_1$ is $u$ (defined in Claim~1 and mentioned in Observation~2).  Such an inbranching 
exists as $C_i$ is strongly connected. Let $L = L_1 \cup L_2 \cup \cdots \cup L_{p-1}$. Note that $V(L)$ and $V(C_p)$
partition $V(D')$. Let $I \subseteq V(L)$ be a set of vertices such that every vertex in $I$ is dominated by at least one vertex 
in $C_p$ (that is, $I \subseteq N^+(C_p) \cap V(L)$). We now define an $I^*$-path as follows.

\2

{\bf Definition ($I^*$-path):} an $I^*$-path is a sequence of vertices $p_0 p_1 p_2 \ldots p_l$, such that the following holds.

\begin{itemize}
\item $p_i \in V(L)$ for all $0 \leq i \leq l$ and $p_l \in I$.
\item $p_i p_{i+1}$ is an $(F-L)$-arc (i.e. $p_i p_{i+1} \in A(F) \setminus A(L)$) 
      or $p_{i+1} p_i$ is a non-$F$-arc (i.e $p_{i+1} p_i \in A(D') \setminus A(F)$) for all $i \in [l-1]$.
\end{itemize}

\begin{figure}
\centering
\begin{tikzpicture}
\draw (0,0) circle [x radius=5cm, y radius=20mm] node {$V(L)$};
\draw (8,0) circle [x radius=15mm, y radius=15mm] node {$C_p$};
\draw (3.5,0) circle [x radius=5mm, y radius=12mm] node {$I$};
\draw [->,color=\greencolor,style=\greenline,style=\thickline] (7,0.5) -- (3.5,0.5);
\draw (3.5,0.65) node {$p_l$};
\draw [->,color=\greencolor,style=\greenline,style=\thickline] (3.5,0.5) -- (2.5,0);
\draw (2.2,0.2) node {$p_{l-1}$};
\draw [->,color=\bluecolor,style=\blueline,style=\thickline] (2.0,-1) -- (2.5,0);
\draw (2.0,-1.15) node {$p_{l-2}$};
\draw [.,dotted,thick] (-2.0,-1) -- (2.0,-1);
\draw (-2.0,-1.15) node {$p_{3}$};
\draw (-2.2,0.2) node {$p_{2}$};
\draw [->,color=\greencolor,style=\greenline,style=\thickline] (-2.0,-1) -- (-2.5,0);
\draw (-3.5,0.65) node {$p_1$};
\draw [->,color=\bluecolor,style=\blueline,style=\thickline] (-3.5,0.5) -- (-2.5,0);
\draw [->,color=\bluecolor,style=\blueline,style=\thickline] (-4,-0.5) -- (-3.5,0.5);
\draw (-4,-0.65) node {$p_0$};

\draw [->,color=\greencolor,style=\greenline,style=\thickline] (7,-2) -- (7.5,-2);
\draw (8.5,-2) node {\footnotesize non-$F$-arcs};
\draw [->,color=\bluecolor,style=\blueline,style=\thickline] (7,-2.3) -- (7.5,-2.3);
\draw (8.5,-2.3) node {\footnotesize $F\setminus L$-arcs};

\end{tikzpicture}
\caption{Example of an $I^*$ path} \label{fig:path}
\end{figure}
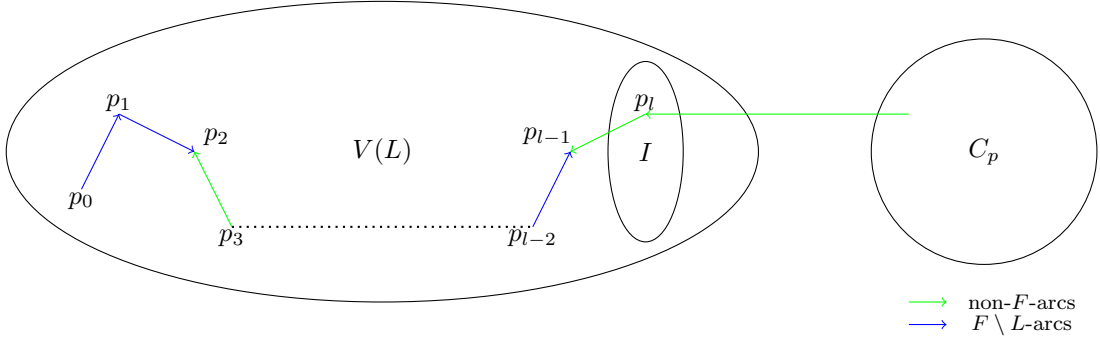

We will now prove the following claims.

\2

{\bf Claim 3:} {\em If $I \subseteq V(L) \cap N^+(C_p)$ and $P=p_0 p_1 p_2 \ldots p_l$ is an $I^*$-path, then the following holds.

\begin{description}
\item[(i):] There is no $w \in V(C_p)$ such that $p_0w,wp_l \in A(D')$.
\item[(ii):] There is no arc $w w^+ \in A(C_p)$ such that $p_0w^+,wp_l \in A(D')$.
\item[(iii):] If $v p_0 \in A(F)$, then there is no $w \in V(C_p)$ such that $vw,wp_l \in A(D')$.
\item[(iv):] If $v p_0 \in A(F)$, then there is no $ww^+ \in A(C_p)$ such that $vw^+,wp_l \in A(D')$.
\end{description}
}

\2

{\em Proof of Claim 3:} Let $I \subseteq V(L) \cap N^+(C_p)$ and let $P=p_0 p_1 p_2 \ldots p_l$ be an $I^*$-path. 
For the sake of contradiction assume that one of the following holds.

\begin{description}
\item[(i):] There is a $w \in V(C_p)$ such that $p_0w,wp_l \in A(D')$. In this case let 
$C^*=p_0p_1\ldots p_lwp_0$.
\item[(ii):] There is an arc $w w^+ \in A(C_p)$ such that $p_0w^+,wp_l \in A(D')$. In this case let 
$C^*=p_0p_1\ldots p_lww^+p_0$.
\item[(iii):] There is a $v p_0 \in A(F)$ and a $w \in V(C_p)$ such that $vw,wp_l \in A(D')$. In this case let 
$C^*=p_0p_1\ldots p_lwvp_0$.
\item[(iv):] There is a $v p_0 \in A(F)$  and a $ww^+ \in A(C_p)$ such that $vw^+,wp_l \in A(D')$. In this case let
$C^*=p_0p_1\ldots p_lww^+vp_0$.
\end{description}

Denote $C*$ as we defined above, we note that $v_i=p_i$ when $i \leq l$ and $c_{l'}=c_0$. Also note that
by the construction either $c_i c_{i+1} \in A(F)$ or $c_{i+1} c_i \in A(D') \setminus A(F)$. We now construct a new eulerian 
factor $F'$ as follows. Initially let $F'=F$ and for all $i \in \{0,1,2,\ldots,l'-1\}$,
if $c_{i+1} c_i \in A(D') \setminus A(F)$ then add $c_{i+1} c_i$ to $F'$ and if this is not the case then we must have
$c_i c_{i+1} \in A(F)$ and in this case remove $c_ic_{i+1}$ from $F'$.
By construction we note that $d_{F'}^+(x) = d_{F'}^-(x)$ for all $x \in V(D')$. Furthermore $F'$ contains all arcs
from $C_p$ except possibly one, so all vertices in $V(C_p)$ belong to the same factor in $F'$. 

In case (i) and (ii) 
above all arcs in $A(L)$ also belong to $A(F')$, so if two vertices belonged to the same factor in $F$ then they also
belong to the same factor in $F'$. As we furthermore have added the arc $wp_l$ to $F'$ we have merged two factors of $F$, implying
that there are fewer factors in $F'$ than in $F$, contradicting the minimality of $p$.  This completes the proof for case
(i) and (ii).

Now assume case (iii) or (iv) holds above. Now all arcs of $L$, except the arc $v p_0$, belong to $A(F')$. Assume
that $v \in C_j$. When we remove the arc $v p_0$ from $L_{j}$ we obtain two connected components $R_v$ and $R_{p_0}$,
where $v \in V(R_v)$ and $p_0 \in V(R_{p_0})$. Furthermore note that the root of $C_j$ belongs to $V(R_{p_0})$.
As we merge the vertices of $R_v$ with $C_p$ (due to the arc $v w$ or $v w^+$) $F'$ also contains at most $p$ factors.
Moreover if $u$ is merged with $C_p$ then $j=1$ and $V(R_{p_0})$ has also been merged with $C_p$ contradicting the minimality of $p$.
So, $u$ has not been merged with $C_p$, implying that the factor of $F'$ containing the vertices of $V(C_p) \cup V(R_v)$ has
more vertices than $C_p$ and does not contain all universal vertices, a contradiction to the maximality of $|V(C_p)|$ in $F$.~\smallqed{}

\2

{\bf Claim 4:} {\em There is no $C_p$-in-universal vertex in $V(L)$.
This also implies that every vertex in $N^+(C_p) \cap V(L)$ has a non-neighbour in $C_p$.

This is in fact true even if $D$ is only $(k+1)$-arc-strong (instead of $\frac{6k+1}{5}$-arc-strong).}

\2

{\em Proof of Claim 4:} Assume for the sake of contradiction that there is a 
$C_p$-in-universal vertex in $V(L)$. Let $I$ be all $C_p$-in-universal vertices in $V(L)$ and
let $X \subseteq V(L)$ contain all vertices $p_0$ such that there exists an $I^*$-path starting at $p_0$.

For the sake of contradiction assume that there is an arc $xy \in A(D') \setminus A(F)$ out of $X$ (that is,
$ x \in X$ and $y \not\in X$). If $y \in V(C_p)$ then we get a contradiction to Claim~3~(i) (where $w=y$ and
the end of the $I^*$-path is $C_p$-in-universal implying that $y$ dominates it). So $y \in A(L) \setminus X$.
However in this case there is an $I^*$-path starting in $y$ (starting with $y x \ldots$), contradicting the 
fact that $ y \not\in X$.  Therefore there is no arc in $V(D') \setminus A(F)$ leaving $X$.

This implies that all arcs of $D'$ leaving $X$ belong to $F$. Assume there are $b$ such arcs. As $D'$ is 
strong we note that $b \geq 1$. As $F$ is an eulerian factor there are also $b$ arcs in $A(F)$
entering $X$.  Let the $F$-arcs entering $X$ be $y_1x_1, y_2x_2, \ldots, y_b x_b$.
Note that $y_1,y_2,\ldots ,y_b$ are distinct (as if $y_i=y_j$ then either $y_ix_i$ or $y_jx_j$ would not belong to $L$ and
$y_i$ would belong to $X$, a contradiction). By the construction of $I$ we note that no $y_i$ is $C_p$-in-universal, implying
that there is a vertex $w_i \in V(C_p)$ which does not dominate $y_i$. By Claim~3~(iii) we note that 
$y_i$ and $w_i$ are non-adjacent in $D'$, implying that there are $b$ pairs of non-neighbours between 
$\{y_1,y_2,\ldots,y_b\}$ and $V(C_p)$, with $b\leq k$.

As there are only $b$ arcs in $D'$ leaving $X$, we note that there must be $k+1-b$ pairs of non-neighbours between 
$X$ and $V(D')\setminus X$ (as $D$ is $(k+1)$-arc-strong). However this implies that there are at least
$k+1-b + b = k+1$ arcs in $A'$, a contradiction. This proves the first part of the claim.

Let $w \in N^+(C_p) \cap V(L)$ be arbitrary. If $w$ has no non-neighbour in $C_p$ then by Observation~2 it is 
either $C_p$-in-universal or $C_p$-out-universal. As $w \in N^+(C_p) \cap V(L)$ it is not $C_p$-out-universal 
and by the above it is not $C_p$-in-universal. Therefore $w$ has a non-neighbour in $C_p$.~\smallqed{}

\2

Recall that $q=\frac{k+1}{5}$.

{\bf Claim 5:} {\em Let $I \subseteq N^+(C_p) \cap V(L)$ such that every vertex in $I$ is dominated 
by at least $i^*$ vertices of $C_p$. Furthermore assume that there are $q-a$ arcs from $C_p$ to 
$V(L) \setminus I$, where $a \geq 0$.

Then there exists a set $W_I \subseteq  V(L) \setminus I$, such that $|W_I| \geq a + |I|$ and
every vertex $w \in W_I$ has at least $i^*+1$ vertices in $V(C_p)$ which are not dominated by $w$.}

\2

{\em Proof of Claim 5:} 
Let $I$, $i^*$ and $a$ be defined as in the statement of the Claim.
Let $X$ be the set of vertices that are the starting point of an $I^*$-path and let $Y = V(L) \setminus X$.
There is no non-$F$-arc (that is arc in $A(D')\setminus A(F)$) that goes from $X$ to $Y$ by the definition of $X$ and $I^*$-paths.
By Claim~4 we note that every vertex in $I$ has a non-neighbour in $C_p$ implying that there are at most
$k-|I|$ non-arcs between $Y$ and $C_p$. Let $z$ denote the number of $F$-arcs that enter $Y$ and note that the
total number of arcs entering $Y$ in $D$ is at most $z + (q-a) + (k-|I|)$. As this number has to be at least $k+q$ we 
get that $z \geq a+|I|$.  Therefore (as there are $z$ $F$-arcs entering $Y$) we observe that there are $z$  $F$-arcs 
leaving $Y$, say $\{y_1x_1,y_2x_2,\ldots,y_z x_z\}$, where $z \geq a + |I|$.
By the definition of $X$, these arcs are $L$-arcs, and as that the $L_i$ are in-branchings, we note that $y_1,y_2,\ldots,y_z$ are distinct
vertices and that $W_I = \{y_1,y_2,\ldots,y_z\}$ therefore has size at least $a + |I|$.

Let $w \in W_I$ be arbitrary and let $p_0 p_1 \ldots p_l$ be an $I^*$-path where $w p_0 = y_i x_i$ for some $i$.
As $p_l \in I$ we note that $N^-(p_l) \cap V(C_p)$ has size at least $i^*$ and as $p_l$ is not $C_p$-in-universal (by 
Claim~4) there is an arc $vv' \in A(C_p)$ entering $N^-(p_l) \cap V(C_p)$. 
By Claim~3~(iii) and (iv) we note that there are no arcs from $w$ to $(N^-(p_l) \cap V(C_p)) \cup \{v\}$, completing the
proof of Claim~5.~\smallqed{}

\begin{figure}
\centering
\begin{tikzpicture}
\draw (0,4) circle [x radius=3cm, y radius=15mm] node {$X$};
\draw (0,0) circle [x radius=3cm, y radius=15mm] node {$Y$};
\draw (6,2) circle [x radius=15mm, y radius=3cm] node {$C_p$};
\draw (1.5,4) circle [x radius=5mm, y radius=10mm] node {$I$};
\draw (0,0.75) circle [x radius=1.5cm, y radius=4mm] node {$W_I$};
\draw [->,color=\greencolor,style=\greenline,style=\thickline] (5,1.5) -- (2,0.75);
\draw [->,color=\greencolor,style=\greenline,style=\thickline] (5.05,1.2) -- (2.2,0.45);
\draw [->,color=\greencolor,style=\greenline,style=\thickline] (5.1,0.9) -- (2.3,0.15);
\draw [->,color=\greencolor,style=\greenline,style=\thickline] (5.2,0.6) -- (2.35,-0.15);
\draw (4,-0.25) node {\footnotesize $\leq q-a$ arcs};
\draw [->,color=\redcolor,style=\redline,style=\thickline] (0,0.95) -- (0,2.9);
\draw [->,color=\redcolor,style=\redline,style=\thickline] (0.3,0.95) -- (0.3,2.9);
\draw [->,color=\redcolor,style=\redline,style=\thickline] (0.6,0.95) -- (0.6,2.9);
\draw [->,color=\redcolor,style=\redline,style=\thickline,style=\redline,style=\thickline] (-0.3,0.95) -- (-0.3,2.9);
\draw [->,color=\redcolor,style=\redline,style=\thickline] (-0.6,0.95) -- (-0.6,2.9);
\draw (1.7,2) node {\footnotesize $\geq |I|+a$ arcs};
\draw [<-,color=\bluecolor,style=\blueline,style=\thickline] (-2.80,-0.3) arc (232:119:2.9cm);
\draw [<-,color=\bluecolor,style=\blueline,style=\thickline] (-2.75,0) arc (230:122:2.6cm);
\draw [<-,color=\redcolor,style=\redline,style=\thickline] (-2.65,0.3) arc (228:125:2.3cm);
\draw [<-,color=\bluecolor,style=\blueline,style=\thickline] (-2.50,0.6) arc (226:130:2.0cm);
\draw (-5,2) node {\footnotesize $\geq |I|+a$ arcs};
\draw [->,color=\greencolor,style=\greenline,style=\thickline] (5,2.5) -- (1.6,3.53);
\draw [->,color=\greencolor,style=\greenline,style=\thickline] (5.05,2.8) -- (1.6,3.55);
\draw [->,color=\greencolor,style=\greenline,style=\thickline] (5.1,3.1) -- (1.6,3.57);
\draw [->,color=\greencolor,style=\greenline,style=\thickline] (5.2,3.4) -- (1.6,3.59);
\draw [->,color=\greencolor,style=\greenline,style=\thickline] (5.3,3.7) -- (1.7,4.33);
\draw [->,color=\greencolor,style=\greenline,style=\thickline] (5.4,4.0) -- (1.7,4.35);
\draw [->,color=\greencolor,style=\greenline,style=\thickline] (5.5,4.3) -- (1.7,4.37);
\draw [->,color=\greencolor,style=\greenline,style=\thickline] (5.6,4.6) -- (1.7,4.39);
\draw (3.9,4.8) node {\footnotesize $\geq i^*\times |I|$ arcs};

\draw [->,color=\greencolor,style=\greenline,style=\thickline] (7,-1) -- (7.5,-1);
\draw (8.5,-1) node {\footnotesize non-$F$-arcs};
\draw [->,color=\bluecolor,style=\blueline,style=\thickline] (7,-1.3) -- (7.5,-1.3);
\draw (8.5,-1.3) node {\footnotesize $F\setminus L$-arcs};
\draw [->,color=\redcolor,style=\redline,style=\thickline] (7,-1.6) -- (7.5,-1.6);
\draw (8.5,-1.6) node {\footnotesize $L$-arcs};

\end{tikzpicture}
\caption{Representation of $W_I$} \label{fig:claim5}
\end{figure}
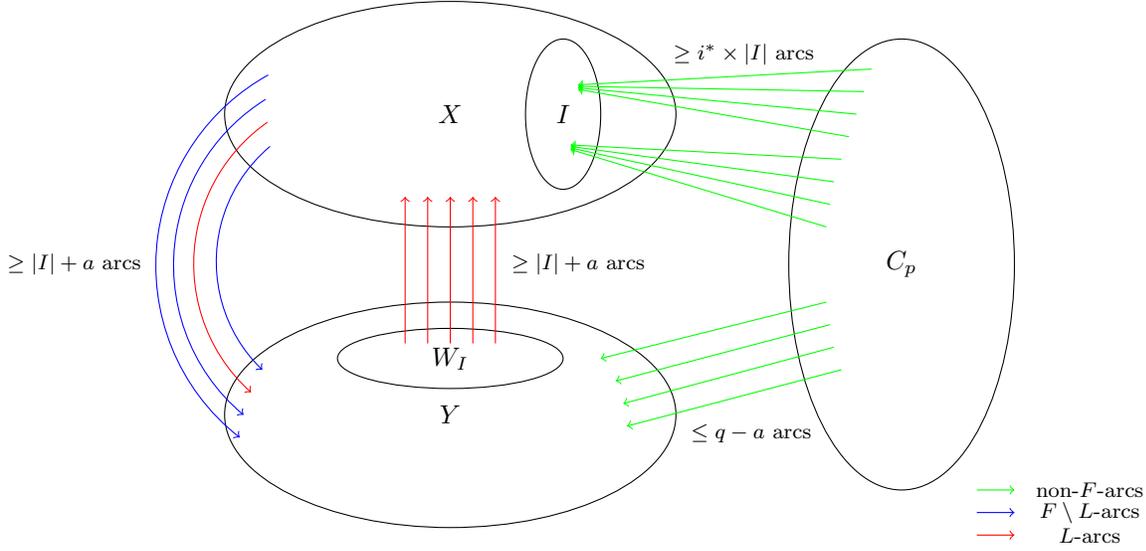

\2

Recall that $q=\frac{k+1}{5}$.

{\bf Claim 6:} {\em There are at least $5q$ non-arcs between $V(L)$ and $V(C_q)$. 

This contradicts $|A'| \leq k$ and thereby completes the proof of the theorem.}

\2

{\em Proof of Claim 6:} Let $\{x_1,x_2,\ldots,x_r\} = N^+(C_p) \cap V(L)$, such that the following holds in $D'$.

\[
d_{C_p}^-(x_1) \leq d_{C_p}^-(x_2) \leq d_{C_p}^-(x_3) \leq \cdots \leq d_{C_p}^-(x_r).
\]

As $D'$ is $q$-arc-strong, we note that there are at least $q$ arcs into $V(L)$ in $D'$.
Therefore $\sum_{i=1}^{r} d_{C_p}^-(x_i) \geq q$.
Define $\theta$ such that it is as large as possible and $\sum_{i=1}^{\theta-1} d_{C_p}^-(x_i) < q$.
Note that $\theta \leq r$.  Let $I_j = \{x_j, x_{j+1}, \ldots, x_r\}$ for $j \in [\theta]$.
Let $s_i = d_{C_p}^-(x_i)$ for all $i \in [\theta]$ and let $t_{\theta} = q - \sum_{i=1}^{\theta-1} s_i$.
Note that $t_{\theta} + \sum_{i=1}^{\theta-1} s_i = q$ and $\sum_{i=1}^{\theta} s_i = q + (s_{\theta}-t_{\theta})$, where $t_{\theta} \leq s_{\theta}$.

We now define a function $f$ as follows.  Initially let $f(v)=0$ for all vertices $v$ in $V(L)$.  
We now modify $f$ as follows.

\begin{itemize}
\item For all $v \in W_{I_1}$ ($I_1$ is defined above and $W_I$ is defined in Claim~5), we let $f(v) = s_1+1$
and for all $v \in I_1$ let $f(v)=1$.
\item For all $v \in W_{I_j}$, for $j=2,3,\ldots,\theta$, increase $f(v)$ by $s_j - s_{j-1}$.
\end{itemize}

We will now show the following two subclaims, where Subclaim~6.2 will complete the proof of Claim~6.

\2

{\bf Subclaim 6.1:} {\em Every $v \in V(L)$ has at least $f(v)$ non-neighbours in $C_p$.}

\2

{\em Proof of Subclaim 6.1:} Let $v \in V(L)$ be arbitrary.  Let $W^* = \cup_{i=1}^{\theta} W_{I_i}$. 
First consider the case when $v \not\in W^*$. If $v \in I_1$, then $f(v)=1$ and by Claim~4  
there is a non-neighbour of $v$ in $C_p$ and therefore Subclaim~6.1 holds in this case.  So we may assume that 
$v \not\in I_1$, in which case $f(v)=0$ and again Subclaim~6.1 holds. We may therefore assume that
$v \in W^*$.

Let $m$ be the minimum value such that $v \in W_{I_m}$ and let $M$ be the maximum value
such that $v \in W_{I_M}$. 
First assume that there are no arcs from $C_p$ to $v$. In this case, by Claim~5, $v$ has
at least $s_M + 1$ non-neighbours in $C_p$. As $f(v) \leq (s_M - s_{M-1})+(s_{M-1}-s_{M-2}) + \cdots
(s_2-s_1) + s_1+1 = s_M+1$, we are done in this case.  So we may assume that there are $s^*$ arcs from
$C_p$ to $v$, where $s^* \geq 1$. As $v \in W_{I_m}$ we note that $v \not\in I_m$ which implies that
$s^* \leq s_m$. Therefore the following holds.

\[
f(v) \leq (s_M - s_{M-1})+(s_{M-1}-s_{M-2}) + \cdots + (s_m - s_{m-1}) + 1 = s_M - s_m + 1 \leq s_M + 1 - s^* 
\]

By Claim~5 we note that there are at least $s_M + 1$ vertices in $C_p$ which are not dominated by $v$, which 
(by the definition of $s^*$) implies that $v$ has at least $s_M + 1 - s^*$ non-neighbours in $C_p$. This
completes the proof of Subclaim~6.1.~\smallqed{}


\2

{\bf Subclaim 6.2:} {\em There are at least $5q$ non-arcs between $V(L)$ and $V(C_q)$}

\2

{\em Proof of Subclaim 6.2:} By the definition of $f$ we note that the following holds.

\[
\begin{array}{rcl} \eqSp{}
\sum_{v \in V(L)} f(v) & = & |I_1| + |W_{I_1}|(s_1+1) + |W_{I_2}|(s_2-s_1) + \cdots + |W_{I_{\theta}}|(s_{\theta}-s_{\theta -1})  \\ 
%
& = & |I_1| + |W_{I_\theta}|(s_{\theta}+1) + (|W_{I_{\theta -1}}|-|W_{I_\theta}|)(s_{\theta-1}+1) + \\
& & \hspace{0.0cm}  (|W_{I_{\theta - 2}}|-|W_{I_\theta - 1}| )(s_{\theta-2}+1)  
+  \cdots + (|W_{I_1}| -|W_{I_2}|)(s_1+1)  \hspace{0.6cm} (*)  \\
 \end{array}
\]

By Claim~5 we note that $|W_{I_j}| \geq (q-\sum_{i=1}^{j-1} s_i)+|I_j|$ for all $j \in [\theta]$.
As $s_1 \leq s_2 \leq \cdots \leq s_{\theta}$ we note that (*) above is minimum when 
$|W_{I_j}| = (q-\sum_{i=1}^{j-1} s_i)+|I_j|$ for all $j \in [\theta]$, which implies that the following holds (as
$\sum_{i=1}^{\theta-1} s_i = q -t_{\theta}$).

\[
\begin{array}{rcl} 
\sum_{v \in V(L)} f(v) & \geq &  |I_1| + (t_{\theta}+|I_{\theta}|)(s_{\theta}+1) 
+ (s_{\theta-1}+|I_{\theta-1}|-|I_{\theta}|)(s_{\theta-1}+1)  \\ \eqSp{}
& & \hspace{0.5cm} + (s_{\theta-2}+|I_{\theta-2}|-|I_{\theta-1}|)(s_{\theta-2}+1)  +  \cdots  + (s_{1}+|I_{1}|-|I_2|)(s_{1}+1) \\
 & \geq & |I_1| + (t_{\theta}+1)(s_{\theta}+1) 
+ (s_{\theta-1}+1)^2
+ (s_{\theta-2}+1)^2
+  \cdots 
+ (s_{1}+1)^2 \\
\end{array}
\]

As $s_{\theta} \geq t_{\theta}$ and $|I_1| \geq \theta$ the above implies the following.

\[
\sum_{v \in V(L)} f(v)  \geq  (t_{\theta}+1)^2 + 1 + \sum_{i=1}^{\theta -1} ((s_i+1)^2 + 1)
\]

Note that for every positive integer $x$ we have $0 \leq (x-1)(x-2) = (x+1)^2+1- 5x$,
which implies that $5x \leq (x+1)^2+1$.  Using the fact that $t_{\theta} + \sum_{i=1}^{\theta -1} s_i = q$
we now obtain the following.

\[
\sum_{v \in V(L)} f(v)  \geq  5t_{\theta} + \sum_{i=1}^{\theta -1} 5s_i = 5q
\]

By Subclaim~6.1, we now note that Subclaim~6.2 holds.~\smallqed{}

\end{proof}

\bibliography{refs}

\end{document}